\newcounter{theor}
\newtheorem{theorem}{Theorem}[theor]

\newtheorem{corollary}{Corollary}[theor]

\newenvironment{proof}[1][Proof]{\begin{trivlist}
\item[\hskip \labelsep {\bfseries #1}]}{\end{trivlist}}
\newenvironment{definition}[1][Definition]{\begin{trivlist}
\item[\hskip \labelsep {\bfseries #1}]}{\end{trivlist}}

\documentclass[12pt]{article}

\begin{document}
\centerline{\bf \large THE SZEG\"O KERNEL ON AN ORBIFOLD CIRCLE
BUNDLE}
\bigskip
\bigskip
\centerline{Jian Song}

\centerline{Department of Mathematics} \centerline{Columbia
University, New York, NY 10027}
\bigskip
\bigskip

\section{\bf  Introduction}
\label{1}
\bigskip

The analysis of holomorphic sections of high powers $L^N$ of
holomorphic ample line bundles $L\rightarrow M$ over compact
K\"ahler manifolds has been widely applied in complex geometry and
mathematical physics. Any polarized K\"ahler metric $g$ with
respect to the ample line bundle $L$ corresponds to the Ricci
curvature of a hermitian metric $h$ on $L$. Any orthonormal basis
$\{S^N_0,...,S^N_{d_N}\}$ of $H^0(M, L^N)$ induces a holomorphic
embedding $\Phi_N$ of M into $CP^{d_N}$. We call the pullback of
the rescaled Fubini-Study metric $\frac{1}{N}\Phi_N^* g_{FS}$ the
Bergman metric with respect to $L^N$. Tian\cite{T1} applied
H\"ormander's $L^2$-estimate to produce peak sections and proves
the $C^2$ convergence of the Bergman metrics. Zelditch\cite{Z}
later generalized Tian's theorem by applying Boutet de
Monvel-Sj\"ostrand\cite{BS} parametrix for the Szeg\"o kernel.
Namely

\begin{theorem}\textnormal{(Zelditch\cite{Z})} Let M be a compact complex manifold of
dimension n and let $(L,h)\rightarrow M$ be a positive Hermitian
holomorphic line bundle. Let g be the K\"ahler metric on M
corresponding to the K\"ahler form $\omega _{g}=Ric(h)$. For each
positive integer $N$, $h$ induces a Hermitian metric $h_{N}$ on
$L^{N}.$ Let $\{S_{0}^{N},S_{1}^{N},...,S_{d_{N}}^{N}\}$ be any
orthonormal basis of $H^{0}(M,L^{N})$, $d_{N}+1=\dim
H^{0}(M,L^{N}),$ with respect to the inner product:
$$
(s_{1},s_{2})_{h_{N}}=\int_{M}h_{N}(s_{1}(x),s_{2}(x))dV_{g},
$$
where $dV_{g}=\frac{1}{n!}\omega _{g}^{n}$ is the volume form of
$g$. Then there is a complete asymptotic expansion:

$$
\sum\limits_{i=0}^{d_{N}}||S_{i}^{N}(z)||_{h_{N}}^{2}\sim
a_{0}(z)N^{n}+a_{1}(z)N^{n-1}+a_{2}(z)N^{n-2}+...
$$
for some smooth coefficients $a_{j}(z)$\ with $a_{0}=1$. More
precisely, for any k:
$$
\left| \left|
\sum\limits_{i=0}^{d_{N}}||S_{i}^{N}(z)||_{h_{N}}^{2}- \sum_{0\leq
j<R}a_{j}(z)N^{n-j}\right| \right| _{C^{k}}\leq C_{R,k}N^{n-R}
$$
where $C_{R,k}$\ depends on $R,k$\ and the manifold $M$.
\end{theorem}

In \cite{L}, Lu shows that each coefficient $a_j(z)$ is a
polynomial of the curvature and its covariant derivatives and
gives a method to compute them explicitly. In particular $a_1(z)$
is the scalar curvature with respect to $g$, which together with
the asymptotic expansion helps Donaldson\cite{D} prove that a
metric of constant scalar curvature on a polarized K\"ahler
manifold is the limit of balanced metrics.

 Since orbifolds arise as degeneration limits of non-singular
K\"ahler manifolds, the property of such limits are crucial to the
understanding of the notion of K-stability conjectured by Tian to
be equivalent to the existence of metrics of constant scalar
curvature. Unfortunately the asymptotic expansions of the Bergman
metrics fail near the singularities in the case of orbifolds. In
this paper we generalize Zelditch's theorem to orbifolds of finite
isolated singularities.

\begin{theorem}

Suppose $M$ is a compact K\"ahler orbifold of $\dim\geq 2$ with
only finite isolated singularities $\{z_i\}_{i=1}^{m}$ and let
$(L,h)\rightarrow M$ be a positive holomorphic orbifold line
bundle. Let g be the orbifold K\"ahler metric on M corresponding
to the K\"ahler form $\omega_g=Ric(h)$. For each N, h induces a
hermitian metric $h_N$ on $L^N$. Let $\{S^N_0,...,S^N_{d_N}\}$ be
any orthonormal basis of $H^0(M, L^N)$, where $d_N+1=\dim H^0(M,
L^N)$, with respect to the inner product
$$<s_1,s_2>_{h_N}=\int_M(s_1(z),s_2(z))h^NdV_g.$$
Let $\{\delta_{z_i}(z)\}_{i=1}^{m}$ be the corresponding
distributions. Then there exists an asymptotic expansion
\begin{eqnarray*}
\sum_{i=0}^{d_N}\|S^N_i(z)\|^2_{h_N}&\sim&
a_0N^n+a_1(z)N^{n-1}+...+a_n(z)+\sum_{i=1}^{m}b(i)\delta_{z_i}(z)\\
&&+a_{n+1}(z)N^{-1}+a_{n+2}(z)N^{-2}+...
\end{eqnarray*}
in the sense that
$$\sum_{i=0}^{d_N}\|S^N_i(z)\|^2_{h_N}-
(a_0N^n+a_1(z)N^{n-1}+...+a_n(z))$$ weakly converges to
$\sum_{i=1}^{m}b(i)\delta_{z_i}(z)$ as $N\rightarrow\infty$.
Furthermore there exist constants $\delta>0$ and $C_{R,k}$ such
that
$$\|\sum_{i=0}^{d_N}\|S^N_i(z)\|^2_{h_N}-
\sum_{0\leq j< R}a_j(z)N^{n-j})\|_{C^k}\leq
C_{R,k}(N^{n-R}+N^{n+k/2}e^{-\delta Nr^2})$$ where $r$ is the
smallest geodesic distance from $z$ to the singularities. In
particular, $a_0=1$ and $a_1(z)$ is the scalar curvature of the
orbifold $(M,g)$ and $b(i)=\frac{1}{|G_i|}\sum_{1\neq g\in
G_i}\frac{1}{\det(I-g|T_{z_i})}$, where $G_i$ is the structure
group of $z_i$ for $i=1,...,m.$
\end{theorem}

We can define the embedding $\Phi_N$: $M\rightarrow
\mathbf{CP}^{d_N}$ by sending $z\in M$ to
$$\Phi_N(z)=[S_0^N(z),..,S_{d_N}^N(z)]\in \mathbf{CP}^{d_N}$$
for $N$ large enough and let $\omega_{FS}$ be the Fubini-Study
metric on $\mathbf{CP}^{d_N}$.

\begin{theorem}
Suppose $M$ is a compact K\"ahler orbifold of $\dim\geq 2$ with
only finite isolated singularities $\{z_i\}_{i=1}^{m}$ and let
$(L,h)\rightarrow M$ be a positive holomorphic orbifold line
bundle. Let $g$ be the K\"ahler metric on M corresponding to the
K\"ahler form $\omega_g=Ric(h)$. For any smooth plurisubharmonic
function $\phi$ with
$\frac{\sqrt{-1}}{2\pi}\partial\overline{\partial}\phi+\omega\geq
0$, we denote $\tilde{h}$ by $he^{-\phi}$ and $\omega_{\tilde{g}}$
by $Ric(\tilde{h})$. Let $\{\tilde{S}^N_0, ... ,
\tilde{S}^N_{d_N}\}$ be any orthonormal basis of $H^0(M, L^N)$
with respect to the inner product
$$<\tilde{s}_1,\tilde{s}_2>_{h_N}=\int_M(\tilde{s}_1(z), \tilde{s}_2(z))\tilde{h}_N \tilde dV_{\tilde{g}}.$$
Then
$$\|\phi-\frac{1}{N}\log(\sum_{i=0}^{d_N}\|\tilde{S}^N_i(z)\|^2_{h_N})\|_{C^0(M)}\rightarrow
0.$$ Furthermore, if we assume that for each $z_i$ its structure
group $G_i$ is abelian, then
$$\|\phi-\frac{1}{N}\log(\sum_{i=0}^{d_N}\|\tilde{S}^N_i(z)\|^2_{h_N})\|_{C^{1,\alpha}(M)}\rightarrow
0,$$ for any positive $\alpha<1$. There also exist $ \epsilon>0$
and $N_0>0$ such that for all $N>N_0$,
$$\inf_{z\in M}\sum_{i=0}^{d_N}\|S^N_i(z)\|^2_{h_N}\geq \epsilon
N^n,$$ where $\{S^N_i\}_{i=0}^{d_N}$ is an orthonormal basis of
$H^0(M, L^N)$ in Theorem 1.2.

\end{theorem}

\begin{corollary}
With the same assumption in as Theorem 2.5, there exist constants
$C_k>0$ and $N_0>0$ such that for any $N>N_0$, we have
$$||\frac{1}{N}\Phi_N^*\omega_{FS}-\omega||_{C^k}\leq C_k (\frac{1}{N}+
N^{k/2}e^{-\delta Nr^2}).$$

\end{corollary}

We conjecture that Theorem 1.3 should be true even without the
assumption that the structure groups be abelian. Tian proved that
any sequence of K\"ahler-Einstein surfaces with positive first
Chern class converges to a K\"ahler-Einstein orbifold and that the
singular points must be rational double points or of cyclic types.
If our conjecture is true then there would exist a uniform
constant $\epsilon>0$ such that for any K\"ahler-Einstein surface
$(M, g)$ with $Ric(g)=g$ we would have
$$\inf_{x\in M}\sum_{i=0}^{d_N}\|S^N_i(z)\|^2_{g}(x)\geq \epsilon N^n.$$

\bigskip
\noindent {\bf Acknowledgements}. The author deeply thanks his
advisor, Professor D.H. Phong for his constant encouragement and
help. He also thanks Professor Zhiqin Lu, Professor Zelditch and
Ben Weinkove for their suggestion on this work. This paper is part
of the author's future Ph.D. thesis in Math Department of Columbia
University.

\section{\bf Orbifolds and orbifold vector bundles}
 \label{2}
\setcounter{theor}{2}
\newtheorem{theorem2}{Theorem}[theor]
\newtheorem{lemma2}{Lemma}[theor]
\newtheorem{corollary2}{Corollary}[theor]
\newtheorem{claim2}{Claim}[theor]
\newtheorem{proposition2}{Proposition}[theor]
\bigskip

\noindent We recall the definition of orbifolds which were
introduced by Satake as $V$-manifolds \cite{B, R, Sa}.

\begin{definition}
An \emph{orbifold structure} on a Hausdorff separable topological
space $X$ is given by an open cover $\bf U$ of X satisfying the
following conditions:
\begin{enumerate}
\item Each $U\in \bf U$ has a local uniformization $\{\tilde{U}, G,
\pi\}$ where $\tilde{U}$ is a connected open neighborhood of the
origin in $\mathbf{C}^n$ and $G$ is a finite group acting smoothly
on $\tilde{U}$ such that $U=\tilde{U} /G$ with $\pi$ as the
projection map. Let $\textrm{ker} G_U$ the subgroup of $G_U$
acting trivially on $U$.

\item If $V \subset U$, then there is a collection of injections
$\{\tilde{V}, G_V, \pi_V\}\rightarrow \{\tilde{U}, G_U, \pi_U\}$.
Namely, the inclusion $i: V \rightarrow U$ can be lifted to
$\tilde{i}: \tilde{V} \rightarrow \tilde{U}$ and an injective
homomorphism $i_{\#}: G_V \rightarrow G_U$ such that $i_{\#}$ is
an isomorphism from $ker G_V$ to $ker G_U$ and $\tilde{i}$ is
$i_{\#}$-equivariant.

\item For any point $x\in U_1\cap U_2$, where $U_1, U_2 \in \bf U$,
there is a $U_3\in \bf U$ such that $x\in U_3\subset U_1\cap U_2$.
\end{enumerate}

\end{definition}

\begin{definition}
An \emph{orbifold bundle} $ B$ over an orbifold $X$  with group
$\Gamma$ and fiber $F$ consists of the following data:
\begin{enumerate}
\item For each local uniformization $\{\tilde{U}, G_U, \pi_U\}$
there is a bundle $B_U$ over $\tilde{U}$ with group $\Gamma$ and
fiber $F$ together with an anti-isomorphism $h_U$ of $G_U$ into a
group of bundle maps of $B_U$ onto itself such that if $b$ lies in
the fiber over $\tilde{x} \in \tilde{U}$, then $h_U(g)b$ lies in
the fiber over $g^{-1} \tilde{x}$ for $g\in G_U$.

\item For any $\tilde{i}: \tilde{V} \rightarrow \tilde{U}$ there
is an induced $i_{\#}$-equivariant bundle map $i^*: B_U
\rightarrow B_V$.

\item If $\tilde{i_1}: \tilde{V} \rightarrow \tilde{U}$ and $\tilde{i_2}: \tilde{U} \rightarrow \tilde{W}$
 then $(i_1 i_2)^*=i_2^* i_1^*$.

 \end{enumerate}
\end{definition}

     Suppose that $X=M/G$ and that $E \rightarrow M$ is a
$G$-equivariant bundle then $E/G \rightarrow X$ is an orbifold
vector bundle.

     The tangent bundle $T_X$ of $X$ is defined by taking for $T_U$
the tangent bundle over $\tilde{U}$, for $h_U(g)$ the inverse of
the mapping of tangent vectors induced by g and $i^*$ the inverse
of the mapping of tangent vectors induced by $i$. If $g$ is a
metric on $T_X$ then for each $\{\tilde{U}, G, \pi\}$, $g_U$ is a
$G$-invariant metric for $\tilde(U)$. We can also define the
cotangent bundle $T_X^*$ and $A^p(X)$ the bundle of differential
p-forms over $X$ in the same manner.

\begin{definition}
Let $\{U_i\}_{i\in I}$ be a locally finite covering of $X$ by open
sets $U_i$ such that $\{\tilde{U_i}, G_{U_i}, \pi_{U_i}\}\in \bf
U$. By a smooth partition of unity for $\{U_i \}_{i \in I}$ we
mean a collection of smooth functions $\{\psi_i\}$ such that
$\textrm{supp}(\psi_i)\subset \tilde{U_i}$ and for each $x\in X$
$\sum_{i\in I} \psi_i(x)=1$.
\end{definition}

   It is easy to show the existence of such partition of unity by
shrinking each $\tilde{U_i}$ a little so that we have a locally
finite covering $\{V_i\}_{i \in I}$ with $V_i\subset U_i$. Then we
can choose for each $i$ a smooth function $u_i$ on $X$ such that
$u_i=1$ on $\tilde{V}_i$ and $u_i=0$ outside $\tilde{U}_i$. Then
we can put $\psi_i=\frac{u_i}{\sum u_i}$.

\noindent {\bf Example 1}

 Let $X$ be the quotient of $\mathbf{CP}^1$ by a cyclic group of order $n$
defined by

$$[Z_0, Z_1]\sim [Z_0 e^{\frac{2k\pi}{n}}, Z_1]$$
for $k=0, 1, 2,...,n-1$. Then it is a "football" which has two
isolated quotient singularities $[0,1]$ and $[1,0]$ with the
cyclic structure group $\mu_{n}$.

\noindent {\bf Example 2}
\begin{definition} Let $d_0,..., d_n$ be $n+1$
positive integers. The weighted projective space $P_{d_0,...,d_n}$
is a toric variety defined by
$$P_{d_0,...,d_n}=\{z \in \mathbf{C}^{n+1}-\{0\} ~| z \sim \lambda z,
\lambda\in \mathbf{C}^* \}$$ where $\mathbf{C}^*$ acts by
$$\lambda (Z_0, ...,Z_n)=(\lambda^{d_0}Z_0,...,\lambda^{d_n}Z_n).$$
\end{definition}
\noindent As for the case of projective spaces, we let
$U_i=\{Z_i\neq 0\}$, then
$$U_i=\{(\frac{Z_0}{Z_i^{d_0/d_i}},...,
\frac{Z_i}{Z_i},...,\frac{Z_n}{Z_i^{d_n/d_i}})\}$$ which is
exactly $\mathbf{C}^n/\mu_{d_i}$ if $d=gcd(d_0,...,d_n)=1.$

   We have the following properties of weighted projective space:
\begin{enumerate}
\item The above $C^*$-action is free if and only if $d_i=d_j$, for
all $i, j=0,...n$.
\item Let $d=gcd(d_0,...,d_n)$ be the greatest common divisor of $d_0,...,d_n$,
then $P_{d_0,...,d_n}$ is homeomorphic to $P_{d_0/d,...,d_n/d}$.
\item Weighted projective spaces are orbifolds which have
singularities with cyclic structure groups acting diagonally. In
particular if $(d_i,d_j)=1$ for all $i\neq j$, $i,j=0,...,n$, then
$P_{d_0,...,d_n}$ has only isolated singularities.
\end{enumerate}

    Now we will state the Riemann-Roch-Kawasaki Theorem which enables
us to determine the coefficients of the currents in the expansion
in Theorem 1.2.

    For each local uniformization $\{\tilde{U},G_U\}$ and each $g\in
G_U$, we consider $\tilde{U}^g$ as a complex manifold on which the
centralizer $Z_{G_U}(g)$ acts. For $V\subset U$, the open
embedding $i:\tilde{V}\rightarrow \tilde{U}$ defines a natural
open embedding $\tilde{V}^h/Z_{G_V}(h)\rightarrow
\tilde{U}^g/Z_{G_U}(g)$ of analytic spaces, where $g=i^{\#}(h)$.
We patch all the $\tilde{U}^g/Z_{G_U}(g)$ together by such
identification which gives a disjoint union of complex orbifolds
of various dimensions:
$$X\amalg\tilde{\Sigma}X=\cup_{\{\tilde{U},G_U\},g\in
G_U}\tilde{U}^g/Z_{G_U}(g).$$ We have a canonical map
$\tilde{\Sigma}X\rightarrow X$ covered locally by the inclusion
$\tilde{U}^g\subset \tilde{U}$. For each $x\in X$ we can choose a
local uniformization $\{\tilde{U}_x, G_x\}$ such that $x\in
\tilde{U}_x$ is a fixed point of $G_x$. $G_x$ is unique up to
isomorphism. Then the number of pieces of $\tilde{\Sigma}X$ is
equal to the number of the conjugacy classes of $G_x$ other than
the identity class.

    Let $\tilde{\Sigma}X_1$,..., $\tilde{\Sigma}X_k$ be all the
connected components of $\tilde{\Sigma}X$. We define $m_i$ for
each $\tilde{\Sigma}X_i$ by
$$m_i=|ker[Z_{G_U}(g)\rightarrow Aut(\tilde{U}^g)|.$$
Let $\sum_{g\in G_U}\mathcal{L}^g(U;E_U)$ be the equivariant Todd
form on $X\amalg\tilde{\Sigma}X$ which represents a cohomology
class $\mathcal{L}(X;E)+\mathcal{L}^{\Sigma}(X;E)$ in
$H^*(X\amalg\tilde{\Sigma}X; C)$. Then we have the following
Riemann-Roch-Kawasaki theorem.

\begin{theorem2}\textnormal{\cite{K}}
Let X be a compact complex orbifold and let $E\rightarrow X$ be a
holomorphic orbifold vector bundle. Then we have
$$\chi (X;\mathcal{O}_X(E))=<\mathcal{L}(X:E),
[X]>+\sum\frac{1}{m_i}<\mathcal{L}^{\Sigma}(X;E),[\tilde{\Sigma}X_i]>.$$
\end{theorem2}
 In particular, if $X$ only has isolated singularities
$\{x_j\}_{j=1,..m}\}$, we have $$\chi
(X;\mathcal{O}_X(E))=<\mathcal{L}(X:E),
[X]>+\sum_{i=1}^{m}\sum_{1\neq g\in
G_{x_i}}\frac{1}{\textrm{det}(1-g|T_{x_i})}.$$

Notice that since $g|T_{x_j}$ is orthogonal and has no eigenvalue
of $1$, it follows that $\det(1-g|T_{x_j})>0.$

\bigskip
\bigskip


\section{\bf The $\overline{\partial}_b$-equation}
\label{3} \setcounter{theor}{3}
\newtheorem{theorem3}{Theorem}[theor]
\newtheorem{lemma3}{Lemma}[theor]
\newtheorem{corollary3}{Corollary}[theor]
\newtheorem{claim3}{Claim}[theor]
\newtheorem{proposition3}{Proposition}[theor]
\bigskip

\noindent In this section we will establish the
$\overline{\partial}_b$-equation for orbifolds and obtain
subelliptic estimates which gives the Hodge decomposition for
$\overline{\partial}_b$ operator. We will essentially follow
Folland and Kohn \cite{FK}.

\begin{definition}
Let $X$ be a compact, orientable real orbifold of dimension
$2n-1$. A partially complex structure on $X$ is an
$(n-1)$-dimensional orbifold subbundle $S$ of $CTX$ such that
\begin{enumerate}
\item $S\cap \overline{S}=\{0\}$,
\item if $L$, $L'$ are local sections of $S$ then so is $[L, L']$.
\end{enumerate}
\end{definition}

\begin{definition}
If $X$ is partially complex, we define the orbifold vector bundle
$B^{p,q}$$(0\leq p, q\leq n-1)$ by $B^{p,q}=\Lambda ^p S^*\otimes
\Lambda^q \overline{S}^*$, which we can identify with an orbifold
subbundle of $\Lambda^{p+q}CTX^*$. We denote by
$\mathcal{B}^{p,q}$ the space of smooth sections of $B^{p,q}$, and
we define $\overline{\partial}_b : \mathcal{B}^{p,q}\rightarrow
\mathcal{B}^{p,q+1}$ as follows:
\begin{enumerate}
\item If $\phi\in \mathcal{B}^{p,0}$, then $\overline{\partial}_b
\phi$ is defined by
$$<\overline{\partial}_b \phi,(L_1\Lambda...\Lambda L_p)\otimes V>=V<\phi,
L_1\Lambda...\Lambda L_p>$$ for all sections $L_1, ..., L_p$ of
$S$ and $V$ of $\overline{S}$.

\item If $\phi \in \mathcal{B}^{p,q}$,
\begin{eqnarray*}
&&(q+1)<\overline{\partial}_b \phi,(L_1\Lambda...\Lambda
L_p)\otimes(V_1\Lambda...\Lambda V_{q+1})>\\
&=&\sum_{j=1}^{q+1}(-1)^{j+1}V_j<\phi,(L_1\Lambda...\Lambda
L_p)\otimes(V_1\Lambda...\hat{V}_j...\Lambda
V_{q+1})>\\
&&+\sum_{i<j}(-1)^{i+j}<\phi,(L_1\Lambda...\Lambda
L_p)\otimes([V_i, V_j]\Lambda V_1\Lambda
...\hat{V}_i...\hat{V}_j...\Lambda V_{q+1})>.
\end{eqnarray*}

\end{enumerate}
\end{definition}
\noindent $\overline{\partial}_b \phi$ is well-defined since
$\overline{\partial}_b$ commutes with $G$.
\begin{eqnarray*}
&&(q+1)<g^*\overline{\partial}_b \phi,(L_1\Lambda...\Lambda
L_p)\otimes(V_1\Lambda...\Lambda V_{q+1})>|_{\tilde{x}}\\
&=&(q+1)<\overline{\partial}_b \phi,g_*(L_1\Lambda...\Lambda
L_p)\otimes g_*(V_1\Lambda...\Lambda V_{q+1})>|_{g \tilde{x}}\\
&=&\sum_{j=1}^{q+1}(-1)^{j+1}V_j<g^*\phi,(L_1\Lambda...\Lambda
L_p)\otimes(V_1\Lambda...\hat{V_j}...\Lambda
V_{q+1})>_{\tilde{x}}\\
&&+\sum_{i<j}(-1)^{i+j}<g^*\phi,(L_1\Lambda...\Lambda
L_p)\otimes([V_i, V_j]\Lambda V_1\Lambda
...\hat{V_i}...\hat{V_j}...\Lambda V_{q+1})>|_{\tilde{x}}\\
&=&(q+1)<\overline{\partial}_b g^*\phi,(L_1\Lambda...\Lambda
L_p)\otimes(V_1\Lambda...\Lambda V_{q+1})>|_{\tilde{x}}.
\end{eqnarray*}

   Let $L_1,..., L_{n-1}$ be a local basis for sections of $S$ over
$\tilde{U}$, so $\overline{L}_1,..., \overline{L}_{n-1}$ is a
local basis for sections of $\overline{S}$. We choose a local
section $N$ of $CTX$ such that $L_1,...,
L_{n-1},\overline{L}_1,..., \overline{L}_{n-1}, N$ span $CTX|_U$
and we may assume that $N$ is purely imaginary. Then the matrix
$(c_{ij})$ defined by
$$[L_i, \overline{L}_j]=\sum a_{ij}^k L_k+\sum
b_{ij}^k \overline{L}_k+c_{ij}N$$ is hermitian and it is called
Levi form.

\begin{proposition3}
The number of non-zero eigenvalues and the absolute value of the
signature of $(c_{ij})$ at each point $\tilde{x}$ are independent
of the choice of $L_1, ..., L_{n-1}, N$.
\end{proposition3}

    We say that $X$ satisfies condition $Y(q)$ if the Levi form have
$\max(q+1, n-q)$ eigenvalues of the same sign or $\min(q+1,n-q)$
pairs of eigenvalues with opposite signs at each point. Notice
that for $n=2$ and $q=1$ the condition $Y(q)$ is never satisfied
and if $X$ is pseudoconvex then $Y(q)$ is satisfied for $n>2$ and
$0<q<n-1$.

We will now define the Sobolev norms on orbifolds.
\begin{definition}
Let $X$ be an $n$-dimensional compact orbifold. Let
$\{U_{\alpha}\}_{\alpha\in \mathcal{A}}$ be a locally finite
covering of $X$ with their uniformization $\{\tilde{U}_{\alpha}\}$
and coordinate mappings $\varphi_{\alpha}:\;
\tilde{U}_{\alpha}\rightarrow R^n$. Let $\{\psi_{\alpha}\}$ be a
partition of unity subordinate to $\{U_{\alpha}\}$. Then for $s\in
R$ and any $k$-form $\phi$ over $X$ we define
$||\phi||_s^2=\sum_{\alpha}||(\psi_{\alpha}\phi)\circ
\varphi_{\alpha}^{-1}||^2_s.$ The norm $||\:||_s$ is not
instrinsic, but it is independent of the choice of local
coordinate charts, partition of unity and the coordinate mappings
up to equivalence.
\end{definition}

    We can choose a hermitian metric on $CTX$ such that $S$,
$\overline{S}$ and $N$ are orthogonal to each other and we can
then assume that $L_1$, $...$, $L_{n-1}$, $\overline{L}_1$, $...$
, $\overline{L}_{n-1}$, $N$ are orthonormal. We can define the
Sobolev spaces $H^{p,q}_s$ for all real $s$ by completing
$\mathcal{B}^{p,q}$ appropriately and define the adjoint operator
$\overline{\partial}^*_b$ and the Laplacian $\Delta
_b=\overline{\partial}_b\overline{\partial}^*_b+\overline{\partial}^*_b\overline{\partial}_b$.
If $\omega_1,...,
\omega_{n-1},\overline{\omega}_1,...,\overline{\omega}_{n-1},\eta$
is the dual basis to $L_1,..., L_{n-1},\overline{L}_1,...,
\overline{L}_{n-1}, N$, we write $\phi\in\mathcal{B}^{p,q}$ as
$\phi=\sum_{IJ}\phi_{IJ}\omega^I\overline{\omega}^J $ and

$\overline{\partial}_b\phi=(-1)^p\sum_{kIJK}\delta^K_{kJ}\overline{L}_k(\phi_{IJ})\omega^I\wedge\overline{\omega}^K$+terms
of order zero,

$\overline{\partial}^*_b\phi=(-1)^{p+1}\sum_{kIHK}\delta^J_{kH}\overline{L}_k(\phi_{IJ})\omega^I\wedge\overline{\omega}^H$+terms
of order zero.

    We define the hermition form $Q_b$ on
$\mathcal{B}^{p,q}$ by
\begin{eqnarray*}
Q_b(\phi,\psi)&=&(\overline{\partial}_b\phi,\overline{\partial}_b\psi)+
(\overline{\partial}^*_b\phi,\overline{\partial}^*_b\psi)+(\phi,\psi)\\
&=&((\Delta_b+I)\phi,\psi).
\end{eqnarray*}

\begin{lemma3}
If $X$ satisfies condition $Y(q)$, then for all $\phi\in
\mathcal{B}^{p,q}$ with support in $U$ we have
$||\phi||^2_{1/2}\leq C Q_b(\phi,\phi).$
\end{lemma3}
\noindent This is a local subelliptic estimate which can be proved
in the same way for non-singular partial complex manifolds with
the condition $Y(Q)$ satisfied. See \cite{FK}.

\begin{theorem3}
If $X$ satisfies condition $Y(q)$, then for all $\phi\in
\mathcal{B}^{p,q}$ we have $||\phi||^2_{1/2}\leq C
Q_b(\phi,\phi).$
\end{theorem3}
\begin{proof}

   We simply apply the partition of unity and
$$||\phi||^2_{1/2}\leq \sum C_1 ||\psi_i \phi||^2_{1/2}\leq C_2\sum
Q_b(\psi_i \phi, \psi_i\phi)\leq C_3Q_b(\phi,\phi).$$
\end{proof}

   We denote the harmonic space by $\mathcal{H}^{p,q}_b=\{\phi\in
\mathcal{B}^{p,q}:\Delta_b \phi=0\}$ which is finite dimensional.
We then have the following Hodge decomposition
$$H^{p,q}_0=\overline{\partial}_b\overline{\partial}^*_b
Dom(\Delta_b)\oplus\overline{\partial}_b\overline{\partial}^*_bDom(\Delta_b)\oplus\mathcal{H}^{p,q}_b.$$
Let $H_b$ be the orthogonal projection on $\mathcal{H}^{p,q}_b$
and $G_b$ be the inverse of $\Delta_b$ on
$(\mathcal{H}^{p,q})^{\bot}$ and zero on $\mathcal{H}^{p,q}_b$. By
the same argument in \cite{FK} we have the following theorem as in
the smooth case.

\begin{theorem3}
Suppose $X$ satisfies condition $Y(q)$ then:
\begin{enumerate}
\item $G_b$ is a compact operator.
\item For any $a\in H^{p,q}_0$, $a=\overline{\partial}_b\overline{\partial}^*_b a
+\overline{\partial}^*_b\overline{\partial}_b a+ H_b a$.
\item $G_bH_b=H_b G_b=0$; $G_b\Delta_b=\Delta_bG_b=I-H_b$ on
$Dom(\Delta_b)$; and if $G_b$ is also defined on
$H^{p,q+1}_0$$(H^{p,q-1}_0)$,
$G_b\overline{\partial}_b=\overline{\partial}_bG_b$ on
$Dom(\overline{\partial}_b)$$(G_b\overline{\partial}^*_b=\overline{\partial}_b^*G_b$
on $Dom(\overline{\partial}^*_b))$.
\item $G_b\mathcal{B}^{p,q}\subset \mathcal{B}^{p,q}$ and $||G_b
a||_s\leq C||a||_{s-1}$ holds uniformly for $a\in
\mathcal{B}^{p,q}$ for each positive integer $s$.
\end{enumerate}
\end{theorem3}

\begin{corollary3}
If $X=\partial M$ is pseudoconvex, then the Szeg\"o projector on
$X$  is given by
$$S=I-\overline{\partial}_b^*G_b\overline{\partial}_b.$$
\end{corollary3}
\begin{proof}
For $n>2$ and $q=1$ $X$ satisfies the condition $Y(1)$ since $X$
is pseudoconvex and we have by the previous theorem that
$$I=\overline{\partial}_b\overline{\partial}^*_bG_b
+\overline{\partial}^*_b\overline{\partial}_bG_b + H_b.$$ So
$\overline{\partial}_bf=\overline{\partial}_b
\overline{\partial}^*_b G_b \overline{\partial}_b f$ where $G_b$
is the Green operator. Thus
$\overline{\partial}_b\overline{\partial}^*_bG_b$ is the
orthogonal projector onto  $Im \overline{\partial}_b$ and
$\overline{\partial}^*_b G_b \overline{\partial}_b$ is the
orthogonal projector onto $(Ker \overline{\partial}_b)^{\bot}$ in
$L^2(X)$. Therefore the Szeg\"o projector $S$ which is the
orthogonal projection from $L^2(X)$ to
$ker(\overline{\partial}_b)\cap L^2(X)$ can be written as
$$S=I-\overline{\partial}^*_b G_b \overline{\partial}_b.$$

\end{proof}

\bigskip
\bigskip
\bigskip


\section{\bf From line bundle to circle bundle}
\label{4} \setcounter{theor}{4}
\newtheorem{theorem4}{Theorem}[theor]
\newtheorem{lemma4}{Lemma}[theor]
\newtheorem{corollary4}{Corollary}[theor]
\newtheorem{claim4}{Claim}[theor]
\newtheorem{proposition4}{Proposition}[theor]
\bigskip

\noindent Let $\mathcal{O}(1)\rightarrow \mathbf{CP}^n$ be the
hyperplane line bundle and let $<,>$ be its natural hermitian
metric. Let $M\in \mathbf{CP}^n$ be a projective manifold and let
$L$ be the restriction of  $M$ and h be the restriction of $<,>$
to $L$. The following lemma is due to Grauert.

\begin{lemma4}
Let D=$\{(z, v)\in L^{\ast}: h(v,v)\leq1\}$. Then D is a strictly
pseudoconvec domain in L.
\end{lemma4}

   Here $L^{\ast}$ is the dual line bundle to $L$. The boundary of
$D$ is a principal $S^1$ bundle $X\rightarrow M$ defined by $
\rho: L^{\ast} \rightarrow R$, $\rho(z,v)=1-|v|^2_z$, where $v\in
L^{\ast}_z$ and $|v|_z$ is its norm in the metric induced by h.
$D=\{\rho>0\}$. We will denote the $S^1$ action by $r_{\theta} x$
and its infinitesimal generator by
$\frac{\partial}{\partial\theta}$ and $\rho$ is $S^1$-invariant.

   Now replace $M$ by an K\"ahler orbifold on which there is an
positive orbforld line bundle $L\rightarrow M$ equipped with the
orbifold hermitian metric $h$.

\begin{lemma4}
Let $D=\{(z, v)\in L^*: h(v,v)\leq1\}$ and $X=\partial D$. Then
both D and X  have an orbifold structure. Furthermore, $X$ is an
orbifold circle bundle over $M$.
\end{lemma4}

   The defining function $\rho(z,v)=1-|v|^2$ is globally defined on
the orbifold $L^*$. Since
$-\frac{\sqrt{-1}}{2\pi}\partial\overline{\partial}\log
h^*=Ric(h^*)<0$ we have
$\frac{\sqrt{-1}}{2\pi}\partial\overline{\partial}h^*>0.$
Therefore on each small enough local uniformization $\rho$ is
convex with respect to a choice of trivialization of $L^*$ and
coordinates on $M$ since $Ric(h)>0$.

\begin{proposition4}
There exists a smooth function $\psi(x,y)$  on each local
uniformization  $\tilde{U}$ such that
\begin{enumerate}
\item $\psi(x,x)=\frac{1}{i}\rho(x)$
\item $d''_x \psi(x,y)$ and $d'_y\psi(x,y)$ vanish on the diagonal
$\{x=y\}$ to infinite order.
\item $\psi(x,y)=-\overline{\psi(y,x)}.$
\end{enumerate}
\end{proposition4}
\noindent By making $\tilde{\psi}(x,y)=\frac{1}{|G|}\sum_{g\in
G}\psi(gx,gy)$ and since $\rho(x)$ is $G$-invariant we can assume
$\psi(.,.)$ is invariant under the diagonal action by $G$.

\begin{proposition4}
There exists a constant $C>0$ such that
$$\textrm{Im}\psi(x,y)\geq C(d(x,X)+d(y,X)+|x-y|^2)+O(|x-y|^3)$$
\end{proposition4}
\begin{proof}
$$\psi(x+h,x+k)\sim \frac{1}{i}\sum
\frac{\partial^{\alpha+\beta}\rho}{\partial
z^{\alpha}\partial\overline{z}^{\beta}}(x)\frac{h^{\alpha}\overline{k}^{\beta}}{\alpha!\beta!}.$$
$$\frac{1}{i}[\psi(x,y)+\psi(y,x)-\psi(x,x)-\psi(y,y)]=L_{\rho}(x-y)+O(|x-y|^3).$$
Also
$$Im\psi(x,y)=\frac{1}{2i}(\psi(x,y)-\overline{\psi(x,y)})=\frac{1}{2i}(\psi(x,y)+\psi(y,x)).$$
\end{proof}

   Notice $\psi$ is only locally defined and in general $\psi$
cannot be globally defined as in the smooth case in \cite{BS} due
to the cancelling of the group action. Also we can always assume
$Im\psi\geq 0$ by shrinking the uniformization a little.

  Denote by $T'D$, $T''D\subset TD\otimes C$ the holomorphic and
antiholomorphic subspaces, and define $d'f=df|_{T'}$ and
$d''f=df|_{T''}$ for $f\in C^{\infty}(D)$. Then $X$ inherits an
orbifold CR structure $CTX=T'\oplus T''\oplus
C(\frac{\partial}{\partial \theta})$. Denote by $T'X$ the space of
holomorphic vector fields on $D$ which are tangent to $X$. They
are given in local coordinates by vector fields $ \sum a_j
\frac{\partial}{\partial \overline{z}_j}$ such that $\sum a_j
\frac{\partial}{\partial\overline{z}_j}\rho=0$. \noindent A local
basis is given by the vector fields
$Z^k_j=\frac{\partial}{\partial\overline{z}_j}-(\frac{\partial
\rho}{\partial\overline{z}_k})^{-1} (\frac{\partial
\rho}{\partial\overline{z}_j})
\frac{\partial}{\partial\overline{z}_k}$ for $j\neq k$.

   The Cauchy-Riemann operator on $X$ is defined by
$$\overline{\partial}_b: C^{\infty}(X) \rightarrow C^{\infty}(X,
(T'')^{\ast})$$ $$\overline{\partial}_b f=df|_{T''}.$$

   It's easy to see $T', T'', \partial \theta$ and
$\overline{\partial}_b$ coincide with $S, \overline{S}, N$ and
$\overline{\partial}_b$ in the previous section.

\begin{lemma4}
$[D_j^k, D_m^k]=0.$
\end{lemma4}
\begin{proof} It can be shown by straightforward calculation.
\end{proof}

\begin{lemma4}
The characteristic cone $\Sigma$ of $\overline{\partial}_b$ is the
real cone of $T^*X$ orthogonal to $T''$ and is generated by
$$ \frac{1}{i}d'\rho|X=-\frac{1}{i}d''\rho|X.$$
\end{lemma4}

\begin{lemma4}
$\sigma([Z, W^*])=L_{\rho}(Z, W)$ where $Z$ and $W$ are two
$C^{\infty}$ local sections of $T''$.
\end{lemma4}
\begin{proof}
\noindent If $Z=\sum a_j\frac{\partial}{\partial\overline{z}_j}$,
$W=\sum b_j\frac{\partial}{\partial\overline{z}_j}$, then
\begin{eqnarray*}
&&\sigma([Z,W^*])(\frac{1}{i}d'\rho)=-\sigma([Z,\overline{W}])(\frac{1}{i}d'\rho)\\
&=&-<[Z,\overline{W}],\frac{1}{i}d'\rho>=-<\sum
a_j\overline{\frac{\partial b_k}{\partial
z_j}}\frac{\partial}{\partial z_k},d'\rho>|_X \\
&=&-\sum a_j\overline{\frac{\partial b_k}{\partial
z_j}}\frac{\partial\rho}{\partial z_k}.
\end{eqnarray*}

\noindent Since $<\overline{W},d'\rho>=<\overline{W}, d\rho>=0$ on
$X$ and $Z$ is tangent to $X$ we have
$$L_Z<W, d'\rho>=a_j(\overline{\frac{\partial b_k}{\partial
z_j}}\frac{\partial\rho}{\partial
z_k}+\overline{b_k}\frac{\partial^2\rho}{\partial\overline{z_j}\partial
z_k})=0.
$$
Therefore $$\sigma([Z,W^*])(\frac{1}{i}d'\rho)=L_{\rho}(Z,W).$$

\end{proof}

   The Hardy space $H^2(X)$ is the space of boundary values of
holomorphic functions on $D$ which are in $L^2(X)$, i.e.
$H^2=\textrm{ker}\overline{\partial}_b\cap L^2(X)$. The $S^1$
action commutes with $\overline{\partial}_b$, hence
$H^2(X)=\oplus_{N} H^2_N(X)$, where $H^2_N(X)=\{f\in H^2(X):
f(\tau_{\theta}x)=e^{iN\theta}f(x)\}$.

  A section of $L$ determines an equivariant
function $\hat{s}$ on $L^{\ast}-\{0\}$ by
$\hat{s}(z,\lambda)=<\lambda, s(z)>$ where $z\in M$ and
$\lambda\in L^{\ast}_z$. Similarly, a section $s_N$ of $L^N$
determines an equivariant function $\hat{s}_N$ on $(L^N)^{\ast}-0$
by $\hat{s}^N(z,\lambda)=<\lambda^N,s_N(z)>$.

\begin{lemma4}
The map $s\rightarrow \hat{s}$ is a unitary equivalence between
$H^0(M, L^N)$ and $H^2(X)$.
\end{lemma4}

   We can generalize the above statements to orbifolds and
holomorphic orbifold line bundles without difficulty.

   We denote by $\Pi : L^2(X) \rightarrow H^2(X)$ and $\Pi_N: L^2(X)
\rightarrow H^2_N(X)$ respectively the orthogonal projections.
Their kernels are defined by
\begin{eqnarray*}
\Pi f(x)=\int_X \Pi(x,y)f(y)d\mu(u)\\
\Pi_N f(x)=\int_X \Pi_N (x,y)f(y)d\mu(u).
\end{eqnarray*}

 Let $\{S^N_i\}$ be an orthonormal basis of
$H^0(M, L^N)$.

\begin{proposition4}
$\| S^N_j(z)\|^2_{h_N}=|\hat{S}^N_i(x)|^2$ for any $x$ with
$\pi(x)=z$.
\end{proposition4}
\begin{proof}
\noindent Let $e_L$ be a local $G$-invariant holomorphic section
$e_L$ of $L$ over a local uniformization $\{\tilde{U},G\}.$ It
induces sections $e^N_L$ of $L^N |_{\tilde{U}}$ and let
$S^N_i(z)=f^N_i(z)e^N_L(z)$ for a holomorphic function $f^N_i$ on
$U$. Then
$$\hat{S}^N_i(z,u)=<u^N,S_i^N(z)>=f^N_i(z)<u^N,
e^N_L(z)>=f^N_i(z)a^{\frac{N}{2}}<u, \frac{e_L}{|e_L|}(z)>^N.$$ So
we have $$\hat{S}^N_i(z,\theta)=f^N_i a(z)^{N/2}e^{iN\theta}.$$
Hence
$|\hat{S}^N_i(z,\theta)|^2=a(z)^N|f^N_i(z)|^2=\|S^N_i(z)\|^2_{h_N}$.
\end{proof}

\begin{proposition4}
$\{\hat{S}^N_i\}$ is an orthonormal basis of $H^2(X)$.
\end{proposition4}
\begin{proof}
\noindent Let $dV_g=\frac{\omega^n_g}{n!}$ be the volume form of
$(M,g)$. Then we have
\begin{eqnarray*}
<S^N_i,S^N_j>&=&\int_M h_N(S^N_i,S^N_j)dV_g\\&=&\int_M
a^N(z)f^N_i(z)\overline{f^N_j(z)}dV_g\\&=&\int_X \hat{S}^N_i
\overline{\hat{S}^N_j}d\mu,
\end{eqnarray*}
where $d\mu=\alpha\wedge d\alpha^n/n!=d\theta\wedge
\pi^{\ast}\omega^n_g$ is the $G_U$-invariant volume form on any
local uniformization $\tilde{U}.$
\end{proof}

\bigskip
\bigskip
\bigskip


\section{\bf The local model}
\label{5} \setcounter{theor}{5}
\newtheorem{theore5}{Theorem}[theor]
\newtheorem{lemma5}{Lemma}[theor]
\newtheorem{corollary5}{Corollary}[theor]
\newtheorem{claim5}{Claim}[theor]
\newtheorem{proposition5}{Proposition}[theor]
\bigskip
\noindent In this and the following section we will follow the
method by L. Boutet de Monvel and J. Sj\"ostrand in \cite{BS} with
a little modification near the singularities to prove a similar
formula for the Szeg\"o kernel for pseudoconvex domain with
quotient singularities.

    Let $(x,y)\in R^n=R^p\times R^q$ and $(\xi, \eta)$ be the dual
variable. Let $\Sigma$ be the cone $\{x=\xi=0\}$. Let $D$ be a
system of pseudo-differential operators
$D_0=\frac{1}{i}(\frac{\partial}{\partial x_j}+x_j|D_y|)$,
$j=1,..., p$.  Let $R$ be a linear continuous operator:
$C_0^{\infty}(R^q)\rightarrow C^{\infty}(R^n)$ defined by
$$Rf(x,y)=(2\pi)^{-q}\int
e^{iy\eta-\frac{1}{2}|x|^2|\eta|}(|\eta|/\pi)^{p/4}\hat{f}(\eta)d\eta.$$

    One has $D_0 R=0$, $R^*R=I$ and
$$I\sim RR^*+L_0D_0.$$ where $L_0\in OPS^{-1, -1}(R^n, \Sigma)$.

    The operator $RR^*$ is defined by the oscillatory integral
$$ RR^*f(x,y)= (2\pi)^{-q}\int \int
e^{i<y-y',\eta>+\frac{1}{2}(|x|^2+|x'|^2)|\eta|}(|\eta|/\pi)^{p/2}f(x',y')dx'dy'd\eta.$$

    The phase function is defined by
$$\phi=<y-y',\eta>+\frac{1}{2}(|x|^2+|x'|^2)|\eta|.$$

    Let $\Sigma^0$ be the cone defined by the complex
equation
 $$\sigma(D_0)=0.$$ Then $\Sigma^0\cap\overline{\Sigma^0}$ is
 the complexified cone of the real cone $\{x=\xi=0\}$. The canonical
 relation $C^+_0$ is the complex cone satisfying
$$\begin{array}{c}
y'=y+i/2(|x|^2+|x'|^2) \\
\xi=ix|\eta|\\
\xi '=-ix'|\eta|\\
\eta '=\eta.
\end{array}
$$ It is easy to see that $C^+_0$ is contained in
$\Sigma^0\times\overline{\Sigma^0}$ and contains
diag$(\Sigma\times\Sigma).$

\begin{proposition5}
The canonical relation $C^+_0$ is unique and satisfies
\begin{enumerate}
\item $C^+_0 \subset \Sigma ^0 \times\overline{\Sigma ^0}$
\item The set of all real points of $C^+$ is exactly the diagonal
of $\Sigma\times\Sigma$
\item $C^+_0$ is positive.
\end{enumerate}
\end{proposition5}

\bigskip
\bigskip
\bigskip


\section{\bf The Szeg\"o kernel on $X$}
\label{6} \setcounter{theor}{6}
\newtheorem{theorem6}{Theorem}[theor]
\newtheorem{lemma6}{Lemma}[theor]
\newtheorem{corollary6}{Corollary}[theor]
\newtheorem{claim6}{Claim}[theor]
\newtheorem{proposition6}{Proposition}[theor]
\bigskip

\noindent
\begin{proposition6}\textnormal{\cite{Bo}}
Let $\{Z_j\}$, $j=1,..,v$ be homogeneous pseudo-differential
operators of degree $m$ and $\Sigma$ the characteristic cone. If
$\{Z_j\}$ satisfy
\begin{enumerate}
\item $\sigma\{[Z_j^*,Z_k]\}$ is positive definite on $\Sigma$
\item $[Z_j, Z_k]\sim \sum A_{jk}^iZ_i$,
\end{enumerate}
then there exists an elliptic Fourier integral operator
transforming the left ideal generated by the $Z_j$ into the ideal
generated by the $\frac{1}{i}(\frac{\partial}{\partial
x_j}+ix_j|D_y|)$, $j=1,...,v$.

\end{proposition6}

   Therefore, for any $\xi\in \Sigma$, there exists a canonical
isomorphism $\Phi$ defined on a neighborhood of $\xi$ and an
elliptic Fourier integral operator $V$ associated with $\Phi$ on a
neighborhood of $\xi$ such that $ \overline{\partial}_b=V^{-1} C
D_0 V$ on a neighborhood of $\xi$, where $C$ is a matrix of
elliptic pseudo-differential operators.

\begin{proposition6}
There exists one canonical relation $C^+$ which is almost analytic
on $T^* \tilde{U}\times T^*\tilde{U}$, unique up to equivalence
and satisfies
\begin{enumerate}
\item $C^+ \subset \Sigma ^0 \times\overline{\Sigma ^0}$.
\item The set of all real points of $C^+$ is exactly the diagonal
of $\Sigma\times\Sigma$.
\item $C^+$ is positive.
\end{enumerate}
\end{proposition6}

\begin{proposition6}
$C^+$ is the canonical relation associated with the phase function
$t\psi(x,y)$ on $\tilde{U}\times \tilde{U}\times R_+$.
\end{proposition6}
\begin{proof} Observe that
\begin{enumerate}
\item $t\psi$ has no critical points
\item $\frac{\partial}{\partial t}t\psi=\psi=0$ implies that the
real critical points is diag$(X)\times R_+$ and on diag$(X) $
$d_x\psi=-d_y\psi=\frac{1}{i}d'\rho|X\neq 0$
\item Im$\psi\geq 0$
\item The set of real points of $C^+$ is exactly diag$\Sigma^+$
as above.
\item Since $d''_x \psi$ and $d'_y\psi$ vanish to infinite order
on diagonal, $C^+$ is contained in
$\Sigma^0\times\overline{\Sigma^0}.$
\end{enumerate}
\end{proof}

\begin{proposition6}
For each point $x$ of $X$ there exists a neighborhood of $x$ with
its local uniformization $\tilde{U}$ and regular operators $S$ and
$L$ such that
\begin{enumerate}
\item $ S\sim S^*\sim S^2$
\item $ \overline{\partial}_b S\sim 0$ and $ I\sim L\overline{\partial}_b+S$,
\end{enumerate}
where $S$ is uniquely determined up to an operator of degree
$-\infty$ and $S\in I^0_c(\tilde{U}^2, C^+)$.

\end{proposition6}

\begin{proof} If $A$ and $B$ are two operators defined on
distribution on $\tilde{U}$. then by $A\sim B$ we mean $R=A-B$ is
an operator of order $-\infty,$ or the kernel distribution of $R$
is $C^{\infty}$ on $\tilde{U}\times\tilde{U}.$

\medskip \noindent {\bf Uniqueness}: Suppose there exist $S'$ and
$L'$ which also satisfy conditions 1 and 2. The assertion of the
theorem is local. Suppose $W$ is an open cone of
$T^{\ast}\tilde{U}-0$, then on $W$
$$S \sim ( S'+ L'\overline{\partial}_b) S \sim S'S$$
and $S' \sim S'^{\ast}$ on $W$. We have on $W$
$$ S' \sim SS' \sim (SS')^{\ast}= S'^{\ast}S^{\ast}\sim S'S\sim S.$$

\noindent {\bf Local existence}

 Let $S_1=V^{-1}RR^*V$, $L_1=V^{-1}L_0C^{-1}V.$ Then we have

$S_1^2=V^{-1}RR^*RR^*V\sim V^{-1}RR^*V\sim S_1$,

$\overline{\partial}_b S_1\sim V^{-1}CD_0 V S_1\sim
V^{-1}CD_0VV^{-1}RR^*V\sim V^{-1}CD_0RR^*V\sim 0$,
\begin{eqnarray*}
&&S_1+L_1\overline{\partial}_b\\
&\sim& V^{-1}RR^*V+V^{-1}L_0C^{-1}VV^{-1}CD_0V\\
&\sim& V^{-1}RR^*V+V^{-1}L_0D_0V\\
&\sim& V^{-1}RR^*V+V^{-1}(I-RR^*)V\sim I.
\end{eqnarray*}

   Since $R^*(V^{-1})^*V^{-1}R$ is an elliptic pseudo-differential
operator, we denote by $B$ its parametrix. Let
$S=V^{-1}RBR^*(V^{-1})^*.$ Then $B\sim B^*$ so $S\sim S^*$.

\noindent
$SS_1=(V^{-1}RBR^*(V^{-1})^*)(V^{-1}RR^*V)=V^{-1}R(BR^*(V^{-1})^*V^{-1}R)R^*V\sim
S_1$

\noindent $S_1S=(V^{-1}RR^*V)(V^{-1}RBR^*(V^{-1})^*)\sim
V^{-1}RBR^*(V^{-1})^*\sim S$

\noindent $\overline{\partial}_bS\sim
\overline{\partial}_bS_1S\sim 0$

\noindent $S^2\sim SS_1S\sim S_1S\sim S.$
 Let $L=(I-S)L_1$, then we have
$$S+L\overline{\partial}_b=S+(I-S)L_1\overline{\partial}_b\sim
S+(I-S)(I-S_1)\sim S+I-S-S_1+SS_1\sim I.$$ \noindent Thus we show
the existence of $S.$

   Let $\{W_{\alpha}\}$ be an open covering of
$T^*\tilde{U} -\{0\}$ and suppose we have $S_{\alpha}$ and
$L_{\alpha}$ satisfying (1) and (2) on each $W_{\alpha}$. By the
same argument as in the proof of uniqueness we can show that on
$W_{\alpha}\cap W_{\beta}$ we have $S_{\alpha} \sim S_{\beta}$.
Using partition of unity as in section 2 we have $Q_{\alpha}$ with
$\sum_{\alpha} Q_{\alpha} \sim I$ and $Q_{\alpha}\sim 0$ outside
$W_{\alpha}$.

   Let $S=\sum Q_{\alpha}S_{\alpha}$ and $L=\sum{\alpha}L_{\alpha}$.
Then on $W_{\beta}$ we have
$$S=\sum Q_{\alpha}S_{\alpha}\sim Q_{\alpha}S_{\beta}\sim
S_{\beta}$$ which gives $S\sim S^*\sim S^2$. Also we have on each
$W_{\beta}$
$$L\overline{\partial}_b=\sum Q_{\alpha}L_{\alpha}\overline{\partial}_b
\sim \sum Q_{\alpha}(I-S_{\alpha})\sim \sum Q_{\alpha}(I-S)\sim
I-S.$$ This completes the proof of the proposition.

\end{proof}

   Such an $S$ admits the following integral representation
$$Sf(x)=\int_{\tilde{U}}\int_0^{\infty}e^{it\psi(x,y)}a(x,y,t)f(y)dtdy,$$
where $f$ has support in $\tilde{U}$ and $a$ is an symbol of
degree $n$:
$$a(x,y,t)\sim \sum_{k=0}^{\infty}t^{n-k}a_k(x,y).$$
The kernel of $S$ can be written as
$$S(x,y)=\int_0^{\infty}e^{it\psi(x,y)}a(x,y,t)dt,$$
and it is smooth off the diagonal. However $S$ is only defined on
$\tilde{U}$ because $Sf(x)$ is not invariant under the action of
the structure group $G_U.$ Since we wish to have $S$ defined on
$U$ instead of its local uniformization, we define

 \noindent $\tilde{S}=\frac{1}{|G_U|}\sum_{g\in G_U}g S g^{-1}$
i.e.

\noindent $\tilde{S} f(x)=\frac{1}{|G_U|}\sum_{g\in G_U} S
\tilde{f}(gx)$ \noindent where $\tilde{f}(x)=f(g^{-1}x).$

\noindent Now $\tilde{S}$ admits the following integral
representation

\begin{eqnarray*}
\tilde{S}f(x)&=&\sum_{g\in
G_U}\int_{\tilde{U}}\int_0^{\infty}e^{it\psi(gx,y)}a(gx,y,t)f(g^{-1}y)dtdy\\
&=&\sum_{g\in
G_U}\int_{\tilde{U}}\int_0^{\infty}e^{it\psi(gx,gy)}a(gx,gy,t)f(y)dtdy\\
&=&\sum_{g\in
G_U}\int_{\tilde{U}}\int_0^{\infty}e^{it\psi(x,y)}a(gx,gy,t)f(y)dtdy.
\end{eqnarray*}
The last equality holds since the volume form $dy$ is
$G_U$-invariant and $\psi$ is invariant under the diagonal action
of $G_U$ . The kernel can be written as
$$\tilde{S}(x,y)=\sum_{g\in G_U}\int_0^{\infty}e^{it\psi(x,y)}a(gx,gy,t)dt
=\int_0^{\infty}e^{it\psi(x,y)}\tilde{a}(x,y,t)dt$$  which is also
smooth off the diagonal. We also define
$\tilde{L}=\frac{1}{|G_U|}\sum g L g^{-1}.$

   Now since $\tilde{S}$ is $G_U$-equivariant, it can be
considered as an operator defined on $U$ instead of on
$\tilde{U}$. The same is true for $L$.

\begin{lemma6} On $\tilde{U}$ we still have

\begin{enumerate}

\item $\tilde{S}\sim \tilde{S}^{\ast} \sim \tilde{S}^2$
\item $\overline{\partial}_b \tilde{S} \sim 0$ and $I \sim
\tilde{L}\overline{\partial}_b + \tilde{S},$
\end{enumerate}
where $A\sim B$ means $R=A-B$ is an operator of order $-\infty$
and is $G_U$-invariant (or equivalently its distribution kernel is
$C^{\infty}$ on $U\times U$).
\end{lemma6}
\begin{proof}

$$\overline{\partial}_{b}\tilde{S}=\frac{1}{|G_U|}\sum_{g\in
G_U}\overline{\partial}_b g S g^{-1} = \frac{1}{|G_U|}\sum g
\overline{\partial}_b S g^{-1} \sim 0$$

$$\tilde{L}\overline{\partial}_b=\frac{1}{|G_U|}\sum g L
\overline{\partial}_b g^{-1}  \sim \frac{1}{|G_U|}\sum g (I-S)
g^{-1}=I-\tilde{S}$$

\begin{eqnarray*}
\tilde{S}^{\ast}&=&\frac{1}{|G|}\sum g S g^{-1}
=\frac{1}{|G_U|}\sum (g S g^{-1})^{\ast}\\
&=&\frac{1}{|G_U|}\sum (g^{-1})^{\ast} S^{\ast} g^{\ast} \sim
\frac{1}{|G_U|}\sum g S g^{-1} \sim \tilde{S}
\end{eqnarray*}

$$ \tilde{S} \sim \tilde{L} \overline{\partial}_b \tilde{S} +
\tilde{S}^2 \sim \tilde{S}^2.$$

\end{proof}

   Now let $\{U_j\}_{j\in J}$ be an open covering of X with their
uniformization $\{U_j, G_j \}$ with partition of unity $Q_j$ and
let $\tilde{S}_j$ and $\tilde{L}_j$ be the corresponding operators
on $\tilde{U_j}$. Let $S_X=\sum Q_j \tilde{S}_j$ and $L_X=\sum Q_j
\tilde{L}_j$.

\begin{definition}
If $A$ and $B$ are two operators defined on distribution on $X$,
then by $A\sim B$ we mean $R=A-B$ is an operator of $-\infty$, or
equivalently the kernel distribution of $R$ is $C^{\infty}$ on
$X\times X$.
\end{definition}

\begin{lemma6}
Given $U_j$ and $U_k$, on the local uniformization
$\tilde{U_j}\cap \tilde{U_k}$ the phase functions $\psi_j$ and
$\psi_k$ of $\tilde{S}_j$ and $\tilde{S}_k$ are equivalent.
Furthermore we have $\tilde{S}_j \sim \tilde{S}_k $.
\end{lemma6}

  Notice that $\tilde{S_j}$ and $\tilde{S}_k$ have the same
canonical relation on $T^* (\tilde{U_j} \cap \tilde{U_k})$ so they
can be composed with each other. And the claim can be easily shown
by the same argument in the proof of the uniqueness in Proposition
6.4.

   On each $U_k$ we have $S_X=\sum Q_j \tilde{S}_j\sim \sum
Q_j\tilde{S}_k\sim {\sum Q_j} \tilde{S}_k\sim \tilde{S}_k$ and
hence
 $$S_X^2\sim S_X^*\sim S_X$$
 $$\overline{\partial}_b S_X\sim 0,$$
 and $$L_X\overline{\partial}_b=\sum Q_j
 \tilde{L}_j\overline{\partial}_b\sim Q_j
 (1-\tilde{S}_j)\sim I - S_X.$$

 So we have proved the following theorem:
\begin{theorem6}
There exist regular operators $S_X$ and $L_X$ on $X$ such that
\begin{enumerate}
\item $ S_X \sim S_X^* \sim S_X^2$
\item $ \overline{\partial}_b \sim 0$ and $ I \sim
L_X\overline{\partial}_b + S_X$,
\end{enumerate}
and $S_X$ is uniquely determined up to an operator of degree
$-\infty$.
\end{theorem6}

   Let $S$ denote the Szeg\"o kernel of $X$. Then, by uniqueness of
$S_X$ and the fact that $I=S+Q_b\overline{\partial}_b$ we have
$S\sim S_X$. For each point $(x_0,x_0)\in X\times X$ we can find a
neighborhood $U$ of $x$ with its local uniformization
$\{\tilde{U},G_U\}$ such that the Szeg\"o kernel $S$ has the local
representation

$$ S_X(x,y)=\sum_{g\in G_U}\int_0^{\infty}e^{it\psi(gx,gy)}a(gx,gy,t)dt,$$
which is smooth off the diagonal. However $S_X(x,y)$ is not
defined on $U\times U$. Remember for any distribution $f$ defined
on $X$, it is on a local uniformization $\{\tilde{U},G_U\}$
$G_U$-invariant. So for $f$ supported on $U$ we have
\begin{eqnarray*}
S_X f(x)&=&\sum_{g\in G_U}\int_{\tilde{U}}\int_0^{\infty}
e^{it\psi(gx,gy)}a(gx,gy,t)f(y)dtdy\\
&=&\sum_{g,h\in G_U}\int_{\tilde{U}}\int_0^{\infty}
e^{it\psi(gx,gy)}a(gx,gy,t)\frac{1}{|G_U|}f(hy)dtdy\\
&=&\sum_{g,h\in G_U}\int_{\tilde{U}}\int_0^{\infty}
e^{it\psi(gx,hy)}a(gx,hy,t)\frac{1}{|G_U|}f(y)dtdy.
\end{eqnarray*}
This enables us to rewrite the Szeg\"o kernel as
$$ \Pi(x,y)=\frac{1}{|G_U|}\sum_{g, h\in
G_U}\int_0^{\infty}e^{it\psi(gx,hy)}a(gx,hy,t)dt.$$And such
$\Pi(x,y)$ is well-defined on $U\times U.$ Although the set of
singularities of $\Pi(x,y)$ sit off-diagonal on
$\tilde{U}\times\tilde{U}$, $\Pi(x,y)$ is smooth off the diagonal
of $U \times U.$

\begin{corollary6} For each point $(x_0,x_0)\in X\times X$ we can
find a neighborhood $U$ of $x$ with its local uniformization
$\{\tilde{U},G_U\}$ such that there exist smooth functions
$F(x,y)$ and $G(x,y)$ on $\tilde{U}\times\tilde{U}$ such that the
Szego kernel has the following representation
$$S(x,y)=\sum_{g,h\in G_U} (F(gx,
hy)(-i\psi(gx,hy))^{-n}+G(gx,hy)\log(-i\psi(gx,hy))).$$

\end{corollary6}

\bigskip
\bigskip
\bigskip


\section{\bf Proof of the main theorems}
\label{7} \setcounter{theor}{7}
\newtheorem{theorem7}{Theorem}[theor]
\newtheorem{lemma7}{Lemma}[theor]
\newtheorem{corollary7}{Corollary}[theor]
\newtheorem{claim7}{Claim}[theor]
\newtheorem{proposition7}{Proposition}[theor]

\bigskip
\noindent On any uniformization $(\tilde{U},G)$ we choose a local
holomorphic coframe $e^{\ast}_L$ and let $a(z)=|e^{\ast}_L|^2_z$
and $(x,y)=(z,\lambda, w,\mu)$ on $X\times X$, we have
$\rho(z,\lambda)=a(z)|\lambda|^2$ and
$$\psi(z,\lambda,w,\mu)=\frac{1}{i}a(z,w)\lambda\overline{\mu},$$
where $a(z,w)$ is an almost analytic function on
$\tilde{U}\times\tilde{U}$ satisfying $a(z,z)=a(z)$.

   On $X$ we have $a(z)|\lambda|^2=1$, so we can assume that
$\lambda=a(z)^{-\frac{1}{2}} e^{\theta}$, then
$$\psi(z,\theta_1,w,\theta_2)=\frac{1}{i}(\frac{a(z,w)}{i\sqrt{a(z)}\sqrt{a(w)}}e^{i(\theta_1-\theta_2)}-1).$$
The weight space projections $\Pi_N$ are Fourier coefficients of
$\Pi$ and hence can be written as
\begin{eqnarray*}
\Pi_N(x,y)&=&\sum_{g,h\in G}\int_{t=0}^{\infty}\int_{S^1}
e^{-iN\theta-it\psi(r_{\theta}gx,hy)}s(r_{\theta}gx, hy,t)d\theta
dt\\
&=&\sum_{g,h\in G}\int_{t=0}^{\infty}\int_{S^1}
e^{iN(-\theta+t\psi(r_{\theta}gx,hy))}s(r_{\theta}gx, hy,
Nt)d\theta dt.
\end{eqnarray*}

   In particular on the diagonal $x=y$, we have
$$\psi(r_{\theta}gx,hx)=\frac{1}{i}(\frac{a(gz,hz)}{a(z)}e^{i\theta}-1).$$
So the phase
$\psi(t,\theta;gx,hx)=\frac{t}{i}(\frac{a(gz,hz)}{a(z)}e^{i\theta}-1)-\theta$.
If $g=h$, then
$\Psi(t,\theta;x,x)=\frac{t}{i}(e^{i\theta}-1)-\theta.$ We have
$d_t\Psi(t,\theta;x,x)=\frac{1}{i}(e^{i\theta}-1)$ and
$d_{\theta}\Psi(t,\theta;x,x)=te^{i\theta}-1$ thus the critical
set is $\{\theta=0, t=1\}.$ The Hessian $\Psi''$ on the critical
set is equal to $\left(
\begin{array}{cc}
0 & 1 \\
1 & i
\end{array}\right).$ So the phase is nondegenerate and the critical points are independent of $x$
and we can apply the theorem in \cite{H2}.  However if $g\neq h$
there is no critical point except $z=0$ and we cannot apply the
theorem and this makes the asymptotic expansion fail near the
singularities.
\begin{eqnarray*}
\Pi_N(z)&=&\sum_{g,h\in G}\int_{t=0}^{\infty}\int_{S^1}
e^{iN(-\theta+t\psi(r_{\theta}gx,hx))}s(r_{\theta}gx,
hx,Nt)d\theta
dt\\
&=&\sum_{g= h\in G}\int_{t=0}^{\infty}\int_{S^1}
e^{-iN\theta-it\psi(r_{\theta}gx,hx)}s(r_{\theta}gx, hx)d\theta
dt\\
&&+\sum_{g\neq h\in G}\int_{S^1} e^{-iN\theta}
\frac{F_{g,h,\theta}}{(1-\frac{a(gz,hz)}{a(z)}e^{i\theta})^{n+1}}
d\theta\\
&&+\sum_{g\neq h\in G}\int_{S^1} e^{-iN\theta}
E_{g,h,\theta}\log(1-\frac{a(gz,hz)}{a(z)}e^{i\theta})d\theta\\
&=& \Pi^{(1)}_N(z)+\Pi^{(2)}_N(z)+\Pi^{(3)}_N(z).
\end{eqnarray*}
 $\Pi^{(1)}_N$ has a converging expansion similar in \cite{Z}
while $\Pi^{(2)}_N$ might cause difficulties near the singular
points.
\begin{eqnarray*}
\Pi^{(2)}_N(z) &=&\sum_{g\neq h\in G}\int_{S^1} e^{-iN\theta}
\frac{F_{g,h,\theta}}{(1-\frac{a(gz,hz)}{a(z)}e^{i\theta})^{n+1}}d\theta\\
&=&\sum_{g\neq h\in G}\int_{S^1}e^{-iN\theta}F_{g,h,\theta}
(\sum_{k=0}^{\infty}\frac{(k+n)!}{n!k!}(\frac{a(gz,hz)}{a(z)})^ke^{ik\theta})d\theta\\
&=&\sum_{g\neq h\in G}\sum_{k=-\infty}^{N}\frac{(l+n)!}{n!l!}
f^{(k)}_{g,h}(\frac{a(gz,hz)}{a(z)})^{N-k},
\end{eqnarray*}
and
\begin{eqnarray*}
\Pi^{(3)}_N(z) &=&\sum_{g\neq h\in G}\int_{S^1} e^{-iN\theta}
E_{g,h,\theta}\log(1-\frac{a(gz,hz)}{a(z)}e^{i\theta})d\theta\\
&=&\sum_{g\neq h\in G}\int_{S^1}e^{-iN\theta}E_{g,h,\theta}
(\sum_{k=0}^{\infty}\frac{1}{k}(\frac{a(gz,hz)}{a(z)})^ke^{ik\theta})d\theta\\
&=&\sum_{g\neq h\in G}\sum_{k+l=N}\frac{1}{l}
e^{(k)}_{g,h}(\frac{a(gz,hz)}{a(z)})^l
\end{eqnarray*}
Here we assume $F_{g,h,\theta}=\sum_k f^{(k)}_{g,h} e^{ik\theta}$
and $E_{g,h,\theta}=\sum_k e^{(k)}_{g,h} e^{ik\theta}$.

\begin{lemma7}
There exists $\delta>0$ such that for any positive integers $l$
and $s$ there is a constant $C_{l,s}$ such that
\begin{enumerate}
\item $||\Pi_N^{(2)}||_{C^s}\leq C_{s,l}
(\frac{1}{N^l}+ N^{n+s/2}e^{-\delta Nr^2}).$
\item $||\Pi_N^{(3)}||_{C^s}\leq C_{s,l}
(\frac{1}{N^l}+ N^{s/2-1}e^{-\delta Nr^2}).$
\end{enumerate}
\end{lemma7}

\begin{proof}
\begin{eqnarray*}
|\Pi_N^{(2)}| &\leq&\sum_{g\neq  h\in G}\sum_{k=-\infty}^N
\frac{(N-k+n)!}{n!(N-k)!}|f^{(k)}_{g,h}|
|\frac{a(gz,hz)}{a(z)}|^{N-k}\\
&\leq& \sum_{g\neq h\in G}\sum_{|k|\geq \sqrt{N}}
\frac{(N-k+n)!}{n!(N-k)!}|f^{(k)}_{g,h}|
|\frac{a(gz,hz)}{a(z)}|^{N-k}\\
&&+\sum_{g\neq h\in G}\sum_{|k|\leq \sqrt{N}}
\frac{(N-k+n)!}{n!(N-k)!}|f^{(k)}_{g,h}|
|\frac{a(gz,hz)}{a(z)}|^{N-k}\\
&\leq& \sum_{g\neq h\in G}\{C_1 \sum_{|k|\geq \sqrt{N}}
(N-k+n)^{n}|f^{(k)}_{g,h}| +C_2\sum_{|k|\leq \sqrt{N}}
N^n|f_{g,h}^{(k)}||\frac{a(gz,hz)}{a(z)}|^{N-k}\}.
\end{eqnarray*}
Since $|f^{(k)}_{g,h}|\leq C_{l}k^l$, we have
\begin{eqnarray*}
|\Pi_N^{(2)}|&\leq& C_l\frac{1}{N^l} +C_2 N^{n}
(\sum_{g\neq h\in G}\sum_{|k|\leq \sqrt{N}}|f_{g,h}^{(k)}|e^{(N-k)\log|\frac{a(gz,hz)}{a(z)}|})\\
&\leq&C_l\frac{1}{N^l}+C_3 N^n (\sum_{g\neq h\in G}\sum_{|k|\leq \sqrt{N}}|f_{g,h}^{(k)}|)e^{-\delta N r^2}\\
&\leq& C_{l,0} (\frac{1}{N^l}+N^ne^{-\delta Nr^2}),
\end{eqnarray*}
and
\begin{eqnarray*}
|\Pi^{(3)}_N|&\leq&\sum_{g\neq  h\in G}\sum_{k=-\infty}^N
\frac{1}{N-k}|e^{(k)}_{g,h}|
|\frac{a(gz,hz)}{a(z)}|^{N-k}\\
&\leq& \sum_{g\neq h\in G}\{ \sum_{|k|\geq \sqrt{N}}
\frac{1}{N-k}|e^{(k)}_{g,h}| +C_4\sum_{|k|\leq \sqrt{N}}
\frac{1}{N}|f_{g,h}|^{(k)}
e^{(N-k)\log|\frac{a(gz,hz)}{a(z)}|}\}\\
&\leq& C_l\frac{1}{N^l} +C_4\frac{1}{N}\sum_{g\neq h\in G}
e^{(N-k)\log(|\frac{a(gz,hz)}{a(z)}|}\\
&\leq&C_l\frac{1}{N^l}+C_4 \frac{1}{N}(\sum_{g\neq h\in G}\sum_{|k|\leq \sqrt{N}}|f_{g,h}^{(k)}|) e^{-\delta N r^2}\\
&\leq& C_{l,0} (\frac{1}{N^l}+\frac{1}{N}e^{-\delta Nr^2}).
\end{eqnarray*}

  We have similar bounds for $|\Pi^{(2)}_N|_{C^k}$ and
$|\Pi^{(3)}_N|_{C^k}$. If $N>\frac{1}{r^2}$where $r$ is the
distance from $x$ to the singular set.  Hence $\Pi^{(2)}_N$ and
$\Pi^{(3)}_N$ converges to $0$ uniformly away from the singular
point. This also proves the first part of Theorem 1.2.
\end{proof}

\begin{lemma7}
Suppose $z_0$ is an isolated singularity on $X$. Then there is a
constant $b(z_0)$ such that for any $C^{\infty}$ test function
$\phi(z)$ supported in a small neighborhood of $z_0$, with $z_0$
its only singularity, we have
$$\lim_{N\rightarrow \infty}<\Pi_N^{(2)}(z)+\Pi_N^{3}(z),
\phi(z)>=\lim_{N\rightarrow \infty}\int_M
(\Pi_N^{(2)}(z)+\Pi_N^{3}(z))\phi(z)\omega^n=b(z_0)\phi(z_0),$$
where $b(z_0)=\sum_{1\neq g\in
G_{z_0}}\frac{1}{\det(I-g|T_{z_0})}$.
\end{lemma7}

\begin{proof}
\begin{eqnarray*}
&&\int_{|z|\leq \epsilon}\Pi_N^{(2)}\phi dV \\
&=&\sum_{g\neq h\in G}\int_{|z|\leq\epsilon}\sum_{|k|\leq
\sqrt{N}} \frac{(N-k+n)!}{n!(N-k)!}f^{(k)}_{g,h}
(\frac{1+gz\overline{hz}}{1+|z|^2})^{N-k}\phi dV\\
&=&\sum_{g\neq h\in G}\int_{|z|\leq\epsilon}\sum_{|k|\leq
\sqrt{N}} \frac{(N-k+n)!}{n!(N-k)!}f^{(k)}_{g,h}
(\frac{1+gz\overline{hz}}{1+|z|^2})^{N-k}\phi(0) dV\\
&&+O(|\sum_{g\neq h\in G}\int_{|z|\leq\epsilon}\sum_{|k|\leq
\sqrt{N}} \frac{(N-k+n)!}{n!(N-k)!}f^{(k)}_{g,h}
(\frac{1+gz\overline{hz}}{1+|z|^2})^{N-k}|z| dV|)\\
&=&\{\sum_{g\neq h\in G}\sum_{|k|\leq
\sqrt{N}}\int_{r=0}^{\epsilon}\int_{S^{2n-1}}\frac{(N-k+n)!}{n!(N-k)!}(1-\frac{r(1-\xi
g\overline{h}\overline{\xi}))}{1+r})^{N-k}
d\xi r^{n-1}dr f^{(k)}_{g,h}(0)\}\phi(0)\\
&&+O(|\sum_{g\neq h\in G}\sum_{|k|\leq
\sqrt{N}}\int_{r=0}^{\epsilon}\int_{S^{2n-1}}\frac{(N-k+n)!}{n!(N-k)!}|1-\frac{r(1-\xi
g\overline{h}\overline{\xi})}{1+r}|^{N-k}d\xi r^ndr|)
\end{eqnarray*}
Let $s=\frac{r}{1+r}$. We have

\begin{eqnarray*}
&&\int_{|z|\leq \epsilon}\Pi_N^{(2)}\phi dV \\
&=&\{\int_{S^{2n-1}}\sum_{g\neq h\in G}\sum_{|k|\leq
\sqrt{N}}\int_{s=0}^{\frac{\epsilon}{1+\epsilon}}\frac{(N-k+n)!}{n!(N-k)!}(1-(1-\xi
g\overline{h}\overline{\xi})s)^{N-k}
s^{2n-1} ds f^{(k)}_{g,h}(0)d\xi \}\phi(0)\\
&&+O(\frac{1}{N})\\
&=&\{\int_{S^{2n-1}}\sum_{g\neq h\in G}\sum_{|k|\leq
\sqrt{N}}\frac{(N-k+n)!}{n!(N-k)!}\frac{n!}{(N-k)...(N-k+n-1)}(1-\xi
g\overline{h}\overline{\xi})^{-(n-1)}\\
&&.\int_{s=0}^{\frac{\epsilon}{1+\epsilon}}(1-(1-\xi
g\overline{h}\overline{\xi})s)^{N-k+n-1}ds f^{(k)}_{g,h}(0)d\xi \}\phi(0)+O(\frac{1}{N})\\
&=&\{\int_{S^{2n-1}}\sum_{g\neq h\in G}\sum_{|k|\leq
\sqrt{N}}(1-\xi g\overline{h}\overline{\xi}))^{n}f^{(k)}_{g,h}(0)d\xi \}\phi(0)+O(\frac{1}{N})\\
&=&\{\sum_{g\neq h\in G}(\int_{S^{2n-1}}\frac{1}{(1-\xi
g\overline{h}\overline{\xi})^n}d\xi)(\sum_{|k|\leq
\sqrt{N}}f^{(k)}_{g,h}(0)) \}\phi(0)+O(\frac{1}{N}).
\end{eqnarray*}
Taking $N\rightarrow \infty$, we have
\begin{eqnarray*}
&&\int_{|z|\leq \epsilon}\Pi_N^{(2)}\phi dV \\
&\rightarrow&(\sum_{g\neq h\in G}\int_{S^{2n-1}}\frac{1}{(1-\xi
g\overline{h}\overline{\xi})^n}d\xi)F((z_0,0),(z_0,0))\phi(z_0)\\
&=&b(z_0)\phi(z_0).
\end{eqnarray*}
The last equation comes from the fact that $F((z_0,0),(z_0,0))=1$.
Also we have $\int_{|z|\leq \epsilon}\Pi_N^{(3)}\phi
dV=O(\frac{1}{N})$ which converges to $0$ as $N\rightarrow
\infty$.
\end{proof}

   We know that $\Sigma_N^{(1)}$ has an asymptotic expansion:
$||\Sigma_N^{(1)}-(a_0 N^n+a_1(z)N^{n-1}+...+a_n(z))||_{C^k}\leq
C_k N^{-1}.$ So $\Pi_N(z)-(a_0 N^n+a_1(z)N^{n-1}+...+a_n(z))$
converges to $\sum_{i=1}^m b(i)\delta_{z_i}(z).$ Therefore we
prove Theorem 1.2.

\begin{corollary}
$\chi_N=\frac{1}{N}\log (\Pi_N) $ converges to $0$ in $C^{0}(M).$
\end{corollary}
\begin{proof}
It suffices to show that $\Pi_N$ is bounded from below by a
uniform positive constant for $N$ large. This can be shown by
constructing equivariant peak sections as in \cite{T2, T3}.
\end{proof}

   Actually we can obtain much stronger result for the special case
where the structure groups are all abelian and from now on we
assume $X$ only has finite isolated singularities with abelian
structure groups.

\begin{lemma7}
If $G$ is a finite abelian subgroup of $U(n)$ then there exists a
uniform constant $C_{G,n}>0$ such that for integer $N>0$ and any
$z\in C^n$ we have
$$\sum_{g\in G}(\frac{1+gz\overline{z}}{1+|z|^2})^N>C_{G,n}.$$

\end{lemma7}
\begin{proof}
Since $G$ is a finite abelian group of isometry of $T_{z_0}$, $G$
can be linearized as a finite abelian subgroup of $U(n)$ and
therefore all elements of $G$ can be diagonalized at the same
time. Any element $g \in G$ can be expressed as
$$g=\left(
\begin{array}{cccc}
e^{\frac{i2p_1(g)\pi}{q_1(g)}} & 0                         & ...& 0 \\
 0                       & e^{\frac{i2p_2(g)\pi}{q_2(g)} }  & ...& 0 \\
 0                       & 0                         & ...& 0\\
 0                       & 0                         & ...& e^{\frac{i2p_n(g)
 \pi}{q_n(g)}}
\end{array}
\right),$$ where $p_i$ and $q_i$ are relatively prime for
$i=1,...,n$.

\noindent Then
\begin{eqnarray*}
&&\sum_{g\in G}(1+gz\overline{z})^N\\
&=&\sum_{g\in
G}\sum_{\alpha_0+\alpha_1+\alpha_n=N}\frac{N!}{\alpha_0!\alpha_1!...\alpha_n!}
|z_1|^{2\alpha_1}|z_2|^{2\alpha_2}...|z_n|^{2\alpha_n}e^{\frac{i2\alpha_1p_1(g)\pi}{q_1(g)}+\frac{i2\alpha_2p_2(g)\pi}{q_2(g)}
    +\frac{i2\alpha_np_n(g)\pi}{q_n(g)}}\\
&=&\sum_{\alpha_0+\alpha_1+\alpha_n=N}\frac{N!}{\alpha_0!\alpha_1!...\alpha_n!}
|z_1|^{2\alpha_1}|z_2|^{2\alpha_2}...|z_n|^{2\alpha_n}[\sum_{g\in
G}e^{\frac{i2\alpha_1p_1(g)\pi}{q_1(g)}+\frac{i2\alpha_2p_2(g)\pi}{q_2(g)}
    +\frac{i2\alpha_np_n(g)\pi}{q_n(g)}}]\\
\end{eqnarray*}
   Fix $\alpha=(\alpha_0, \alpha_1, ... , \alpha_n)$, we can
construct the following group homomorphism
\begin{eqnarray*}
&&\alpha: G  \rightarrow  U(1)\\
&&\alpha(g)=e^{\frac{i2\alpha_1p_1(g)\pi}{q_1(g)}+\frac{i2\alpha_2p_2(g)\pi}{q_2(g)}
    +\frac{i2\alpha_np_n(g)\pi}{q_n(g)}}.
\end{eqnarray*}

 It is easy to see that $\sum_{g\in
G}e^{\frac{i2\alpha_1p_1(g)\pi}{q_1(g)}+\frac{i2\alpha_2p_2(g)\pi}{q_2(g)}
    +\frac{i2\alpha_np_n(g)\pi}{q_n(g)}}$ is nonzero only if
    $\alpha$ is a trivial homomorphism and in this case $\sum_{g\in
G}e^{\frac{i2\alpha_1p_1(g)\pi}{q_1(g)}+\frac{i2\alpha_2p_2(g)\pi}{q_2(g)}
    +\frac{i2\alpha_np_n(g)\pi}{q_n(g)}}$  is a positive integer. Therefore
\begin{eqnarray*}
&&\sum_{g\in G}(1+gz\overline{z})^N\\
&\geq&\sum_{g\in G}\sum_{
\begin{array}{c}
\alpha_0+\alpha_1+\alpha_n=N\\
|G||\alpha_k, \textnormal{ for $k=1,...,n$}\\
\end{array}}
\frac{N!}{\alpha_0!\alpha_1!...\alpha_n!}
|z_1|^{2\alpha_1}|z_2|^{2\alpha_2}...|z_n|^{2\alpha_n}\\
&\geq& C_{G,n}(1+|z|^2)^N.
\end{eqnarray*}

\end{proof}

\begin{lemma7}
There exist constants $N_0$, $C$ and $c>0$ such that for $N>N_0$
we have for all $z\in M$ $$ cN^n \leq \Pi_N(z)\leq C N^n.$$
\end{lemma7}

\begin{proof}
For the proof of Theorem 1.3, we have the upper bound. For the
lower bound it suffices to prove it near the singularities.
\begin{eqnarray*}
\Pi_N&=&\sum_{g, h\in G}\sum_{k=-\infty}^N
\frac{(N-k+n)!}{n!(N-k)!}f^{(k)}_{g,h}(\frac{a(gz,hz)}{a(z)})^{N-k}\\
&=&\sum_{g, h\in G}\sum_{|k|\leq \sqrt{N}}
\frac{(N-k+n)!}{n!(N-k)!}f^{(k)}_{g,h}(\frac{a(gz,hz)}{a(z)})^{N-k}+O(N^{n-1})\\
&=&\sum_{g, h\in G}\sum_{|k|\leq \sqrt{N}}
\frac{N^n}{n!}f^{(k)}_{g,h}(\frac{a(gz,hz)}{a(z)})^{N-k}+O(N^{n-1})\\
&=&\sum_{g, h\in G}\sum_{|k|\leq \sqrt{N}}
\frac{N^n}{n!}f^{(k)}_{g,h}(z=0)(\frac{a(gz,hz)}{a(z)})^{N-\sqrt{N}}+O(N^{n-\epsilon'})\\
&=&\frac{N^n}{n!}\sum_{|k|\leq \sqrt{N}}f^{(k)}(z=0)\sum_{g\neq
h\in
G}(\frac{a(gz,hz)}{a(z)})^{N-\sqrt{N}}+O(N^{n-\epsilon'})\\
&\geq&C_5N^n\sum_{g \in
G}(\frac{1+gz\overline{z}}{1+|z|^2})^{N-\sqrt{N}}+O(N^{n-\epsilon'})\geq
cN^n
\end{eqnarray*}
by lemma 7.3 and here $\sum_{|k|\leq \sqrt{N}}f^{(k)}(z=0)$
converges to $F(0,0)$ which is positive.

\end{proof}

\begin{proposition7}
$\chi_N=\frac{1}{N}\log\Pi_N $ converges to 0 in $C^{1,\alpha}(M)$
for any $\alpha<1$.
\end{proposition7}
\begin{proof}
$$\bigtriangledown
\chi_N=\frac{1}{N}\frac{\bigtriangledown\Pi_N^{(1)}+\bigtriangledown\Pi_N^{(2)}+\bigtriangledown\Pi_N^{(3)}}{\Pi_N}.$$
Since we have the above lemma we can do the same calculation in
Lemma 7.1 and it is straightforward to prove the lemma.
\end{proof}

\begin{proposition7}
Let
$\omega(N)=\frac{1}{N}\log(\sum_{i=0}^{d_N}\|S^N_I(z)\|^2_{h_N})$
be the pullback of the scaled Fubini-Study metric then
$$||\omega(N)-\omega||_{C^k}\leq C_k (\frac{1}{N}+(\sum_{l=0}^k
N^{k/2}(Nr^2)^{l/2+1})e^{-\delta Nr^2}).$$ In particular,
$||\omega(N)-\omega||_{C^k}\leq C_k' N^{k/2}.$
\end{proposition7}
\noindent The proof is straightforward by induction.

   Let $P(M,\omega)=\{\phi\in
C^{\infty}(M)\;|\;\omega_{\phi}=\omega+\frac{\sqrt{-1}}{2\pi}\partial\overline{\partial}\phi>0,\sup_{M}\phi=0\}$
be the set of all plurisubharmonic functions on $M$. Let
\begin{eqnarray*}
\tilde{\omega _{\phi}} &=&\omega +\frac{\sqrt{-1}}{2\pi}\partial
\overline{\partial} \phi
>0\ \\
\widetilde{h} &=&he^{-\phi }.
\end{eqnarray*}
Let $\widetilde{h}_{N}$ be the induced Hermitian metric of
$\tilde{h}$ on $L^{N}$, $\{\tilde{S}
_{0}^{m},\tilde{S}_{1,...,}^{N}\tilde{S}_{d_{N}}^{N}\}$ be any
orthonormal basis of $H^{0}(M,L^{N})$($1+d_{N}=\dim
H^{0}(M,L^{N})$) with respect to $\tilde{h},\omega_{\phi}$. Then
we have the following holomorphic approximation theorem which has
been proved by Lu.
\begin{corollary7}
$\phi_N=\frac{1}{N}\log
(\sum_{k=0}^{d_{N}}||\tilde{S}^N_k(z)||^2_{h^N}) $ converges to
$\phi$ in $C^0(M)$ for any $\alpha<1$. If for each singularity its
structure group is cyclic, then the $\phi_N$ converges to $\phi$
in $C^{1,\alpha}(M)$ for any $\alpha<1$.
\end{corollary7}

\bigskip
\bigskip
\bigskip


\section {\bf Examples}\label{8}

\bigskip

 \noindent Let $X$ be the quotient of $\mathbf{CP}^1$ by a cyclic group of
order $n$ defined by
$$[Z_0, Z_1]\sim [Z_0 e^{\frac{i2k\pi}{n}}, Z_1]$$
for $k=0, 1, 2,...,n-1$.

Let $L$ be $O(n)$ the orbifold line bundle over $X$ and
 the Fubini-Study metric on $X$ is defined by
$g_{i\overline{j}}=\frac{\sqrt{-1}}{2\pi}\partial_i\partial_{\overline{j}}\log(|Z_0|^2+|Z_1|^2)$.
Then
$\{\sqrt{\frac{n(nN+1)!}{(nk)!(nN-nk)!}}Z_0^{nk}Z_1^{nN-nk}\}_{k=0}^{N}$
is an orthonormal basis of $H^0(X, O(nN))$ with respect to the
Fubini-Study metric.

   On the patch $U_0=\{Z_0\neq 0\}$ we write
$r=\frac{|Z_1|^2}{|Z_0|^2}$ and we have
\begin{eqnarray*}
S^{(N)}(z)&=&n\sum_{m=0}^{N}(nN+1)C_{nN}^{nm}
\frac{|Z_0|^{2nm}|Z_1|^{2nN-2nm}}{(|Z_0|^2+|Z_1|^2)^{nN}}\\
&=&n\sum_{k=m}^{N}(nN+1) C_{nN}^{nm}\frac{r^{nm}}{(1+r)^{nN}}\\
&=&(nN+1)\sum_{k=0}^{n-1}(\frac{1+re^{\frac{i2k\pi}{n}}}{1+r})^{nN}\\
&=&nN+1+nN
\sum_{k=1}^{n-1}(\frac{1+re^{\frac{i2k\pi}{n}}}{1+r})^{nN}.
\end{eqnarray*}
And
$$|S^{(N)}|([1,0])=|S^{(N)}|^2([0,1])=n(nN+1)$$
$$\lim_{N\rightarrow\infty,[Z_0, Z_1]\neq [0,1],[1,0]}\frac{|S^{(N)}|^2([Z_0,Z_1])}{nN+1}=1.$$

    Also we have
\begin{eqnarray*}
|S^{(N)}(z)-(nN+1)|&=&(nN+1)|\sum_{k=1}^{n-1}(\frac{1+re^{\frac{i2k\pi}{n}}}{1+r})^{nN}|\\
&\leq&n(nN+1)\max|1-\frac{r(1-e^{\frac{i2k\pi}{n}})}{1+r}|^{nN}\\
&\leq&n(nN+1)e^{-\delta Nr}\\
&\leq&C_2\frac{1}{N},
\end{eqnarray*}
if $N\geq \frac{1}{r^2}$.

   By similar calculation we have
$$S^{(N)}\sim
nN+1+\frac{n-1}{2n}([\{[0,1]\}]+[\{[1,0]\}]),$$ as a distribution
on $X$ and if we integrate $S^{(N)}$ over $X$ we obtain the
Gauss-Bonnet theorem for the quotient sphere. It also verifies
Theorem 1.2 since $a(0)=1$, $a(1)=1$ as the scalar curvature and
$$b_{[0,1]}=b_{[1,0]}=\frac{1}{n}\sum_{k=1}^{n-1}\frac{1}{1-e^{\frac{i2k\pi}{n}}}=\frac{n-1}{2n} .$$

\end{document}